\documentclass[aip,
amsmath,amssymb,
reprint,%
numerical,
]{revtex4-1}
\usepackage{graphicx}% Include figure files
\usepackage{dcolumn}% Align table columns on decimal point
\usepackage{bm}% bold math
\usepackage{hyperref}
\usepackage[utf8]{inputenc}
\usepackage[T1]{fontenc}
\usepackage{mathptmx}
\usepackage{etoolbox}
\newcommand{\RNum}[1]{\uppercase\expandafter{\romannumeral #1\relax}}
\usepackage{array,xcolor,colortbl}
\usepackage{multirow}
\makeatletter
\def\@email#1#2{%
	\endgroup
	\patchcmd{\titleblock@produce}
	{\frontmatter@RRAPformat}
	{\frontmatter@RRAPformat{\produce@RRAP{*#1\href{mailto:#2}{#2}}}\frontmatter@RRAPformat}
	{}{}
}%
\makeatother

\begin{document}
	
	\title{Effect of heating or cooling in a suspension of phototactic algae with no-slip boundary conditions}
	% Force line breaks with \\
	
	% Force line breaks with \\
	\author{S. K. Rajput}
	\altaffiliation[Corresponding author: E-mail: ]{shubh.iiitj@gmail.com.}%Lines break automatically or can be forced with \\
	\author{M. K. Panda}%
	%\email{mkpanda@iiitdmj.ac.in}
	%\email{shubh.iiitj@gmail.com}
	\affiliation{$^1$ Department of Mathematics, PDPM Indian Institute of Information Technology Design and Manufacturing, Jabalpur 482005, India.%\\This line break forced with \textbackslash\textbackslash
	}%
	%\homepage{http://www.Second.institution.edu/~Charlie.Author.}
	%\affiliation{%
	%	Second institution and/or address%\\This line break forced% with \\
	%}%
	

	%%%%%%%%%%%%%%%%%%%%%%%%%%%%%%%%%%%%%%%%%%%%%%%%%%%%%%%%%%%%%%%%%%%%%%%%%%%		
	\begin{abstract}
		
		In this study, we investigate the impact of heating or cooling in a suspension experiencing phototactic bioconvection. The suspension is illuminated by collimated irradiation from the top and subjected to heating or cooling from the bottom. The governing equations include the Navier–Stokes equations with the Boussinesq approximation, the diffusion equation for motile microorganisms, and the energy equation for temperature. Employing linear perturbation theory, we analyze the linear stability of the suspension. The findings predict that the suspension undergoes destabilization when heated from below and stabilization when cooled from below. This suggests a sensitive dependence of the system's stability on the thermal conditions, providing valuable insights into the behavior of phototactic bioconvection under different heating or cooling scenarios.
		
	\end{abstract}
	
	%%%%%%%%%%%%%%%%%%%%%%%%%%%%%%%%%%%%%%%%%%%%%%%%%%%%%%%%%%%%%%%%%%%%%%%%%%%	
	
	\maketitle
	
	%%%%%%%%%%%%%%%%%%%%%%%%%%%%%%%%%%%%%%%%%%%%%%%%%%%%%%%%%%%%%%%%%%%%%%%%%%%	
	
	\section{INTRODUCTION}
	
	Bioconvection, a captivating phenomenon, showcases the convective motion observed in fluid containing self-propelled motile microorganisms like algae and bacteria at a macroscopic level ~\cite{20platt1961,21pedley1992,22hill2005,23bees2020,24javadi2020}. These microorganisms exhibit an intriguing tendency to move upwards on average due to their higher density compared to the surrounding medium, typically water. The formation of distinct patterns in bioconvection is closely tied to the behavior of these motile microorganisms. However, pattern formation is not solely reliant on their upswimming or higher density; rather, it is influenced by their response to various environmental stimuli known as "taxes", including gravitaxis, chemotaxis, phototaxis, gyrotaxis and thermotaxis. Phototaxis is a biological phenomenon in which living organisms move in response to light. On the other hand, in thermotaxis, living organisms or cells exhibit directed movement in response to changes in temperature. 
	
	The term "thermal photactic bioconvection" refers to a specific subset of bioconvection phenomena where both light (photo) and temperature (thermal) gradients play crucial roles in orchestrating the motion and behaviors of microorganisms, typically motile algae. It represents the unique interaction between phototaxis (movement in response to light) and thermotaxis (movement in response to temperature gradients) in these microorganisms. This phenomenon is a fascinating field of study because it allows researchers to explore the synergistic effects of these two environmental cues on the movement and distribution of aquatic microorganisms.

	Kuznetsov~\cite{51kuznetsov2005thermo} conducted a study on the onset of thermo-bioconvection in a dilute suspension of oxytactic microorganisms. In this investigation, a shallow fluid layer of the suspension heated from below was found to be less stable than under isothermal conditions.
	Alloui et al.~\cite{52alloui2006stability} explored the effect of heating or cooling from below on the stability of a suspension of motile gravitactic microorganisms in a shallow fluid layer. They observed that thermal effects could either stabilize or destabilize the suspension, leading to a decrease or increase in the wavelength of the bioconvective pattern.
	Nield and Kuznetsov~\cite{53nield2006onset} employed linear stability analysis to study the onset of bioconvection in a horizontal layer of fluid containing a suspension of motile microorganisms, considering heating or cooling from below. Their findings revealed that the stability boundary depended on the values of the Lewis and Prandtl numbers, and oscillatory convection could be the favoured mode of instability when the layer is heated from the bottom.
	Alloui et al.~\cite{54alloui2007numerical} investigated the effect of heating or cooling from below in a square enclosure on the development of gravitactic bioconvection. They presented the influence of thermo-effects on a bifurcation diagram and flow structure, reporting that heating from below destabilizes the suspension while cooling from below stabilizes it.
	Taheri and Bilgen \cite{55taheri2008thermo} studied the effect of heating or cooling from below, at constant temperature and constant heat flux, on the development of gravitactic bioconvection in vertical cylinders with stress-free sidewalls. They reported that the thermal effect considerably modified the pattern formation of gravitactic bioconvection.
	Kuznetsov \cite{56kuznetsov2011bio} developed a theory of bio-thermal convection in a suspension containing two species of microorganisms, gyrotactic and oxytactic, exhibiting different taxes. He identified one traditional Rayleigh number and two bioconvection Rayleigh numbers for gyrotactic and oxytactic microorganisms. Increasing any of the three Rayleigh numbers decreased the stability of the system.
	Saini et al. \cite{57saini2018analysis} investigated bio-thermal convection in a suspension containing gravitactic microorganisms, saturated by a fluid within the framework of linear and nonlinear stability theory. They determined that the bioconvection Rayleigh number destabilizes the onset of biothermal convection, and this effect is more predominant at high microorganism speeds.
	Zhao et al. \cite{57zhao2018linear} applied linear stability analysis to investigate the stability of bioconvection in a suspension of randomly swimming gyrotactic microorganisms, heated from below. They reported that the Lewis number had no effect on the critical value of the thermal Rayleigh number but had a significant influence on the critical bioconvection Rayleigh number.
	
	In the realm of phototactic bioconvection, the pioneering work of Vincent and Hill~\cite{12vincent1996} marked its initial presence. They delved into the impact of collimated irradiation on an absorbing (non-scattering) medium. Subsequently, Ghorai and Hill~\cite{10ghorai2005} explored the behavior of phototactic algal suspension in two dimensions, excluding the consideration of scattering effects.
	The investigation into light scattering, both isotropic and anisotropic, was undertaken by Ghorai $et$ $al$.\cite{7ghorai2010} and Ghorai and Panda\cite{13ghorai2013} with normal collimated irradiation. Panda and Ghorai~\cite{14panda2013} proposed a model for an isotropically scattering medium in two dimensions, yielding results different from those reported by Ghorai and Hill~\cite{10ghorai2005} due to the inclusion of scattering effects.
	Panda and Singh~\cite{11panda2016} explored phototactic bioconvection in two dimensions, confining a non-scattering suspension between rigid side walls. The impact of diffuse irradiation, in combination with collimated irradiation, was investigated by Panda $et$ $al$.\cite{15panda2016} in an isotropic scattering medium and by Panda\cite{8panda2020} in an anisotropic medium.
	Considering that, in natural environments, sunlight strikes the Earth's surface at off-normal angles with oblique irradiation, Panda $et$ $al$.\cite{16panda2022} and Kumar\cite{17kumar2022} investigated the effects of oblique collimated irradiation on non-scattering and isotropic scattering suspensions. Kumar~\cite{40kumar2023,39kumar2023} further explored the onset of phototaxis in the rotating frame, considering the effects of a rigid top surface and scattering when the suspension was illuminated by normal collimated flux.
	In a recent study, Panda and Rajput~\cite{41rajput2023} delved into the effects of diffuse flux in conjunction with oblique collimated flux on isotropic scattering algae suspension. The open literature has witnessed research addressing the interplay of diffuse and collimated radiation in tandem with natural convection.

	While thermal bioconvection has been studied extensively for its inherent biological significance, in the whole literature the thermal effect of the phototactic bioconvection did not account. This research unveils a fascinating twist to this phenomenon by highlighting the profound impact of temperature gradients on the phototactic response of these organisms. Thermal effects on phototactic bioconvection not only intensify the complexities of this natural marvel but also open up new dimensions of understanding the interplay between temperature, light, and biology.
	
	The structure of the article is organized as follows: Initially, the problem is mathematically formulated, followed by the derivation of a fundamental (equilibrium) solution. Subsequently, a small disturbance is introduced to the equilibrium system, and the linear stability problem is acquired through the application of linear perturbation theory, followed by numerical solution methods. The model results are then presented, and finally, the implications and findings of the model are thoroughly discussed.
	
	%%%%%%%%%%%%%%%%%%%%%%%%%%%%%%%%%%%%%%%%%%%%%%%%%%%%%%%%%%%%%%%%%%%%%%%%%%%	
	
	\begin{figure}[!htbp]
		\centering
		\includegraphics[width=8cm]{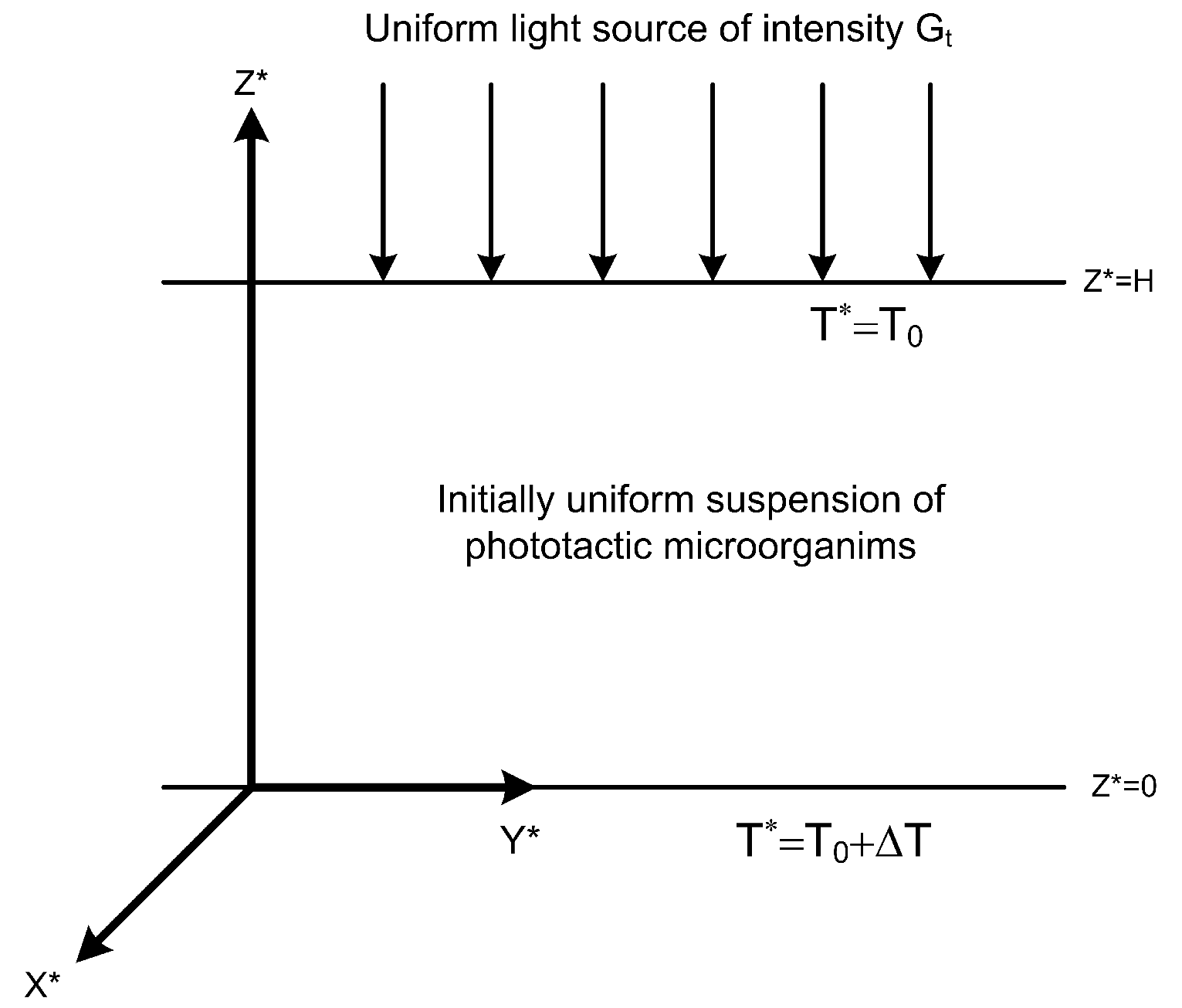}
		\caption{\footnotesize{The geometry of the problem under consideration.}}
		\label{fig1}
	\end{figure}
	
	\section{Phototaxis and thermotaxis in suspensions of algae}
	
	Consider a suspension of phototactic algae of infinite width with horizontal boundaries at $z^*= 0, H$. The lower and upper boundaries are taken as rigid non-slip. Indeed, even if the upper boundary is open to the air, cells often collect at the top surface forming a rigid-like packed layer. This layer is subject to uniform collimated radiation from above and heated and cooled from below. It is assumed that the temperature variation is sufficiently weak so that it does not kill microorganisms and phototactic behavior (cells swimming orientation and speed etc.) of microorganisms remains unaffected due to heating or cooling. The temperature $T^*$, at the lower boundary $(z^*=0)$ is maintained at $T_0+\delta T$ and at the upper boundary $(z^*=H)$ is maintained at $T_0$,
	which are considered to be uniform. 
	
	Let $\boldsymbol{p}$ be the unit vector corresponding to the swimming direction and $<\boldsymbol{p}>$ is the ensemble average of the swimming direction for all the cells in an elemental volume. For many species of micro-organisms, the swimming speed is independent of illumination, position, time and direction and we denote the ensemble-average swimming speed by $W_c$. The average swimming velocity is thus

	\begin{equation*}
	\boldsymbol{W}_c=W_c<\boldsymbol{p}>.
	\end{equation*}
	The mean swimming direction, $<\boldsymbol{p}>$, is given by
	\begin{equation}\label{1}
	<\boldsymbol{p}>=M(G)\hat{\boldsymbol{z}},
	\end{equation}
	where $\hat{\boldsymbol{z}}$, is a unit vector in the vertical $z$-axis and $M(G)$ is the phototaxis function such
	that
	\begin{equation*}
	M(G)=\left\{\begin{array}{ll}\geq 0, & \mbox{ } G(\boldsymbol{x}^*)\leq G_{c}, \\
	< 0, & \mbox{ }G(\boldsymbol{x}^*)>G_{c}.  \end{array}\right. 
	\end{equation*}
	$M(G)$ depends on the total light intensity $G$ reaching the cell and has the generic form shown in Fig. 1. The exact functional form of $M$ will depend on the species of microorganisms. The mean swimming direction is zero when the light intensity is critical or optimal (at $G =G_c$, $M(G)=0$).
	For a suspension of algae subject to unidirectional illumination, the light intensity at any position is a function of cell concentration because microorganisms nearer to the light source shade those that are further away. Then the light intensity, $G(\boldsymbol{x}, t)$, at
	position $\boldsymbol{x}=(x, y, z)$ in the suspension is given by
	
	\begin{equation}\label{2}
	G(\boldsymbol{x}^*,t^*)=G_t\exp\left(-\alpha^*\int_H^{\boldsymbol{ r}}n({\boldsymbol{ r}',t^*})d\boldsymbol{ r}'\right).
	\end{equation}

	Here $\boldsymbol{ r}$ is the vector from the cell at position $\boldsymbol{x}$ to the light source, which has an intensity $G_t$, and $n$ is the concentration of cells at any position $\boldsymbol{x}$. $\alpha$ is the extinction coefficient, which quantifies the strength of absorption and scattering of light in the suspension. Eq.~(3) is derived from the Lambert-Beer law for weak scattering and It was first predicted by Kessler that the Lambert-Beer law could be used as the basis of a
	model for algal phototaxis, although it has been used previously to calculate absorption in algal suspensions by Duysens. We shall assume that light absorption by the pure fluid is negligible (as we have considered the layer to be not deep enough). Also, the algal cell only receives the light that is travelling to it in a direct line from the source, hence the effect of scattering by the algae cells is ignored. Thus this expression for light intensity only describes the absorption by micro-organisms.
	
	%%%%%%%%%%%%%%%%%%%%%%%%%%%%%%%%%%%%%%%%%%%%%%%%%%%%%%%%%%%%%%%%%%%%%%%%%%%	     
	\subsection{The continuum model}
	
	The basic equations governing the flow in the medium for a dilute and incompressible suspension of phototactic microorganisms are based on the continuum model proposed by Vincent and Hill~\cite{12vincent1996}. Here, It is assumed that heating or cooling from below is sufficiently weak, so it does not kill microorganisms and does not affect their phototactic behavior. Each cell has a volume $v$ and density $\rho+\Delta\rho$, where $\rho$ is the density of the fluid in which the cells swim and $\Delta\rho/\rho<<1$. $\boldsymbol{u}^*$ is the average velocity of all the material in a small volume $\delta V$ and $n^*$ is the cell concentration. The governing equations are obtained by using the  Navier–Stokes equations of motion, energy equation, and cell conservation equation, over a representative elementary volume. Under these assumptions, the governing equations for three-dimensional flow in the medium are obtained as
	
	Equation of continuity
	\begin{equation}\label{3}
	\boldsymbol{\nabla}^*\cdot \boldsymbol{u}^*=0,
	\end{equation}
	The momentum balance equation
	\begin{equation}\label{4}
	\rho\left(\frac{D\boldsymbol{u}^*}{D t^*}\right)=-\boldsymbol{\nabla} P_e+\mu{\nabla^*}^2\boldsymbol{u}^*-n^*v g\Delta\rho\hat{\boldsymbol{z}}+\rho g\beta(T^*-T_0)\hat{\boldsymbol{z}},
	\end{equation}
	Here $D/Dt^* = \partial/\partial t^* + \boldsymbol{u}^* \cdot \boldsymbol{\nabla}^*$ is the material derivative, $P_e$ is the excess pressure above hydrostatic, $\hat{\boldsymbol{z}}$ is a unit vector vertically upward, and $\mu$ is the viscosity of the suspension which is assumed to be that of the fluid and $\beta$ is the thermal expansion coefficient.
	
	The cell conservation equation
	\begin{equation}\label{5}
	\frac{\partial n^*}{\partial t^*}=-\boldsymbol{\nabla}^*\cdot \boldsymbol{F}^*,
	\end{equation}
	
	where $\boldsymbol{F}^*$ is the net flux and it can be written as
	\begin{equation*}
	\boldsymbol{F}^*=n^*(\boldsymbol{u}^*+W_c<\boldsymbol{p}>)-\boldsymbol{D}\cdot\boldsymbol{\nabla} n^*.
	\end{equation*}
	The thermal energy equation
	\begin{equation}\label{6}
	\frac{\partial T^*}{\partial t^*}+\boldsymbol{\nabla}^*\cdot(\boldsymbol{u}^*T^*) =\alpha{\boldsymbol{\nabla}^*}^2 T^*.
	\end{equation}
	Two important assumptions have been made in the expression for the cell flux vector. First, the cells are considered to be purely phototactic and thus the effect of viscous torque, which might contribute to the swimming component, is neglected. Second, the diffusion tensor $D$, which should be derived from the swimming velocity autocorrelation function, is assumed to be constant. These assumptions allow us to remove the Fokker-Planck equation from the governing equations and the resulting model is a valid limiting case to consider to understand the complexity of the problem before moving to a more complex detailed model.
	
	Now, the expression for light intensity,
	$G(\boldsymbol{x}^*, t^*)$, at position $\boldsymbol{x}^*= (x^*, y^*, z^*)$ in the suspension
	
	\begin{equation}\label{7}
	G(\boldsymbol{x}^*,t^*)=G_t\exp\left(-\alpha^*\int_{z^*}^H n^*({\boldsymbol{ x}',t^*})d\boldsymbol{ z}'\right),
	\end{equation}    
	
	provided the absorption across the suspension is weak $(0 < \alpha^* n^*H \ll 1)$. The boundary conditions are 
	
	\begin{subequations}   
		\begin{equation}\label{8a}
		\boldsymbol{u}\cdot\hat{\boldsymbol{z}}=\boldsymbol{F}^*\cdot\hat{\boldsymbol{z}}=0\qquad at\quad z^*=0,H,
		\end{equation}
		and for the rigid boundary
		\begin{equation}\label{8b}
		\boldsymbol{u}^*\times\hat{\boldsymbol{z}}^*=0\qquad at\quad z^*=0,H.
		\end{equation}
		For temperature
		\begin{equation}\label{8c}
		T^*=T_0+\Delta T\qquad at\quad z^*=0,
		\end{equation}
		\begin{equation}\label{8d}
		T^*=T_0\qquad at\quad z^*=H.
		\end{equation}
	\end{subequations}
	
	The governing equations are made dimensionless by scaling all lengths on $H$, the depth of the layer, time on a diffusive time scale $H^2/D$, and the bulk fluid velocity on $D/H$. The appropriate scaling for the pressure is $\mu D/H^2$, the cell concentration is scaled on $\bar{n}$ the mean concentration and temperature are scaled by $\Delta T$. The unstarred symbol indicates the scalled governing parameters. In terms of the non-dimensional variables, the bioconvection equations become
	
	\begin{equation}\label{9}
	\boldsymbol{\nabla}\cdot\boldsymbol{u}=0,
	\end{equation}
	\begin{equation}\label{10}
	S_c^{-1}\left(\frac{D\boldsymbol{u}}{D t}\right)=-\nabla P_{e}+\nabla^{2}\boldsymbol{u}-Ra_bn\hat{\boldsymbol{z}}+Ra_TLeT\hat{\boldsymbol{z}},
	\end{equation}
	\begin{equation}\label{11}
	\frac{\partial{n}}{\partial{t}}=-\boldsymbol{\nabla}\cdot\boldsymbol{F},
	\end{equation}
	where
	\begin{equation*}
	\boldsymbol{F}=\boldsymbol{n{\boldsymbol{u}}+nV_{c}<{\boldsymbol{p}}>-{\boldsymbol{\nabla}}n},
	\end{equation*}
	and 
	\begin{equation}\label{12}
	\frac{\partial T}{\partial t}+\boldsymbol{\nabla}\cdot(\boldsymbol{u}T) =Le\boldsymbol{\nabla}^2 T.
	\end{equation}
	Here, $S_c^{-1}=\nu/D$ is the Schmidt number, $V_c=W_cH/D$ is the dimensionless swimming speed, $Ra=\bar{n}v g\Delta{\rho}H^{3}/\nu\rho{D}$ is bio-convective Rayleigh number,  $Ra_T= g\beta\Delta{T}H^{3}/\nu\alpha$ is the thermal Rayleigh number, and $Le=\alpha/D$ is the Lewis number.
	
	The light intensity, $G(\boldsymbol{x}, t)$, at position $\boldsymbol{x} = (x, y, z)$ in dimensionless form becomes
	\begin{equation}\label{13}
	G(\boldsymbol{x},t)=G_t\exp\left(-\kappa\int_z^1 n({\boldsymbol{ x}',t})d\boldsymbol{ x}'\right),
	\end{equation} 
	where, $\kappa = \alpha n\bar{H}$, is the non-dimensional extinction coefficient.

	After non-dimensionalization, the boundary conditions are expressed as
	
	\begin{subequations}
		\begin{equation}\label{14a}
		\boldsymbol{u}\cdot\hat{\boldsymbol{z}}=\boldsymbol{F}\cdot\hat{\boldsymbol{z}}=0\qquad at\quad z=0,1,
		\end{equation}
		and for the rigid boundary
		\begin{equation}\label{14b}
		\boldsymbol{u}\times\hat{\boldsymbol{z}}=0\qquad at\quad z=0,1.
		\end{equation}
		For temperature
		\begin{equation}\label{14c}
		T=1\qquad at\quad z=0,
		\end{equation}
		\begin{equation}\label{14d}
		T=0\qquad at\quad z=1.
		\end{equation}
	\end{subequations}
	
	%%%%%%%%%%%%%%%%%%%%%%%%%%%%%%%%%%%%%%%%%%%%%%%%%%%%%%%%%%%%%%%%%%%%%%%%%%%	
	
	\section{The steady solution}
	
	Equations (\ref{9})–(\ref{12}) possess a static equilibrium solution in which
	\begin{equation}\label{15}
	\boldsymbol{u}=0,~~~n=n_s(z),\quad T=T_s(z)\quad and\quad  <\boldsymbol{p}>=<\boldsymbol{p}_s>=M(G_s)\hat{\boldsymbol{z}}.
	\end{equation}
	
	The concentration $n_s(z)$ satisfies
	
	\begin{equation}\label{16}
	\frac{dn_s}{dz}-V_cM_sn_s=0,
	\end{equation}
	which is supplemented by the cell conservation relation
	\begin{equation*}
	\int_0^1n_s(z)dz=1,
	\end{equation*}
	where $M_s = M(G)$ at $G = G_s$ . Here, the steady light intensity $G_s$ at a height $z$ $(0 \leq z \leq 1)$ is defined as
	
	\begin{equation}\label{17}
	G_s(\boldsymbol{z})=G_t\exp\left(-\kappa\int_z^1 n_s({\boldsymbol{ z}'})d\boldsymbol{ z}'\right),
	\end{equation}    
	
	Let us introduce a new variable $\bar\omega=\int_1^z n_s({\boldsymbol{ z}'})d\boldsymbol{ z}'$. Now, in terms of this new variable $\bar{\omega}$ Equation (19) and (20) changes to
	
	\begin{equation}\label{18}
	\frac{d^2\bar{\omega}}{dz^2}-V_cM_s\frac{d\bar{\omega}}{dz}=0
	\end{equation}
	
	with the boundary condition
	
	\begin{subequations}\label{19}
		
		\begin{equation}\label{19a}
		\bar{\omega}+1=0,~~~~~at ~~z=0,
		\end{equation}
		\begin{equation}\label{19b}
		\bar{\omega}=0,~~~~~at ~~z=1.
		\end{equation}
		
	\end{subequations}
	
	$T_s(z)$ satisfies
	
	\begin{equation}\label{20}
	\frac{d^2T_s}{dz^2}=0
	\end{equation}
	with the boundary conditions (\ref{14c}) and (\ref{14d}) which implies
	\begin{equation}\label{21}
	T_s(z)=1-z
	\end{equation}
	
	Equations (\ref{18}) and (\ref{20}) with boundary conditions (\ref{19}) and (\ref{21}) respectively constitute a boundary value problem.
	Vincent and Hill are considered the most interesting and simplified case when steady light intensity $G_s\approx G_c$ throughout a horizontal layer of the suspension of unit depth, then the phototaxis function $M_s$ changes to a linear function of $G_s$ after linearizing it about $G_c$ using Taylor series and neglecting the higher order terms. Indeed, they have used the fact that one can approximate a curve by its tangent.
	
	In this study, we have considered a phototaxis function $M_s$ such that the graph it will approximately fit the sketch of the phototaxis function obtained from the experiments. Now, the phototaxis function becomes non-linear and its graph covers the whole suspension depth. The phototaxis function in our study is generated by superimposing the following sine functions:
	
	\begin{equation}\label{22}
	M_s=0.8\sin\left[\frac{3\pi}{2}\Xi(G)\right]-0.1\sin\left[\frac{\pi}{2}\Xi(G)\right],
	\end{equation}
	where $\Xi(G)=G\exp[\beta(G-1)]$. However, the functional form of the phototaxis function is not unique. One can generate its functional form by superimposing other trigonometric functions. The value of critical intensity is related to the parameter $\beta$. When $G_c=G_t$, i.e., the critical intensity occurs at the upper surface, the cells accumulate at the top of the domain. This is similar to the gravitactic cells. As $G_c$ decreases, the maximum concentration decreases and the location of the maximum concentration
	shifts towards the midheight of the domain. The maximum concentration is the smallest when the maximum is located in the middle of the domain. As $G_c$ decreases further, the maximum concentration increases and shifts towards the bottom of the domain. In our study, we have taken the value of the light intensity at source as 0.8 and the range of the $\beta$ has been chosen as $-1.1\leq\beta\leq1.1$, for which the value of the critical intensity $G_c$ lies in the range $0.3\leq G_c\leq G_t$ $(=0.8)$. The range of $\beta$ has been chosen such that the plot of the functional form of the phototaxis function shall fit suitably the experimental sketch of the phototaxis function. 
	
	%%%%%%%%%%%%%%%%%%%%%%%%%%%%%%%%%%%%%%%%%%%%%%%%%%%%%%%%%%%%%%%%%%%%%%%%%%%	   
	
	\section{Linear stability of the problem}
	We consider a small perturbation of amplitude $\epsilon (0 < \epsilon << 1)$ to the equilibrium state (\ref{15}) so that
	
	\begin{equation}\label{23}
	\begin{pmatrix}
	\boldsymbol{u}\\n\\T\\<\boldsymbol{p}>
	\end{pmatrix}
	=
	\begin{pmatrix}
	0\\n_s\\T_s\\<\boldsymbol{p}_s>
	\end{pmatrix}
	+\epsilon
	\begin{pmatrix}
	\boldsymbol{u}_1\\n_1\\T_1\\<\boldsymbol{p}_1>
	\end{pmatrix}
	+O(\epsilon^2),
	\end{equation}
	where  $\boldsymbol{u}_1=(u_1,v_1,w_1)$.
	Substituting the perturbed variables into Eqs. (13)–(15) and linearizing about the equilibrium state by collecting $O(\epsilon)$ term, gives
	\begin{equation}\label{24}
	\boldsymbol{\nabla}\cdot \boldsymbol{u}_1=0,
	\end{equation}
	\begin{equation}\label{25}
	S_c^{-1}\left(\frac{\partial \boldsymbol{u_1}}{\partial t}\right)=-\boldsymbol{\nabla} P_{e}+\nabla^{2}\boldsymbol{u}_1-Ra_bn_1\hat{\boldsymbol{z}}+Ra_TLeT_1\hat{\boldsymbol{z}},
	\end{equation}
	\begin{equation}\label{26}
	\frac{\partial{n_1}}{\partial{t}}+V_c\boldsymbol{\nabla}\cdot(<\boldsymbol{p_s}>n_1+<\boldsymbol{p_1}>n_s)+w_1\frac{dn_s}{dz}=\boldsymbol{\nabla}^2n_1,
	\end{equation}
	\begin{equation}\label{27}
	\frac{\partial{T_1}}{\partial t}+w_1\frac{dT_S}{dz}=Le\boldsymbol{\nabla}^2T_1.
	\end{equation}
	If $G=G_s+\epsilon G_1$ , then the perturbed intensity $G_1$ [after expanding $\exp(-\kappa\int_z^1\epsilon n_1 dz')$ and collecting $\mathcal{O}(\epsilon)$ term] is given by
	
	\begin{equation}\label{28}
	G_1(\boldsymbol{z})=G_s\exp\left(-\kappa\int_z^1 n_1({\boldsymbol{ z}'})d\boldsymbol{ z}'\right).
	\end{equation}  
	
	Now the expression $T(G_s+\epsilon G_1)\hat{\boldsymbol{z}}-<\boldsymbol{ p}_s>$ on simplification at $O(\epsilon)$ gives perturbed swimming direction
	\begin{equation}\label{29}
	<\boldsymbol{p}_1>=G_1\frac{dT_s}{dG}\hat{\boldsymbol{z}}.
	\end{equation}
	
	Now putting the value of $\boldsymbol{ p}_1$ from Eq.~(\ref{29}) into Eq.~(\ref{26}) and simplifying we get
	
	\begin{equation}\label{30}
	\frac{\partial n_1}{\partial t}+w_1\frac{dn_s}{dz}+\Gamma_1(z)\int_1^zn_1 dz'+\Gamma_2(z)n_1+\Gamma_3(z)\frac{\partial n_1}{\partial z}=\nabla^2 n_1, 
	\end{equation}
	
	where
	\begin{subequations}
		
		\begin{equation}\label{31a}
		\Gamma_1(z)=\kappa V_c\frac{d}{dz}\left(n_sG_s\frac{dM_s}{dG}\right),
		\end{equation}
		\begin{equation}\label{31b}
		\Gamma_2(z)=2\kappa V_c n_s G_s\frac{dM_s}{dG},
		\end{equation}
		\begin{equation}\label{31c}
		\Gamma_3(z)=V_cM_s.
		\end{equation}
	\end{subequations}

	By elimination of $P_e$ and the horizontal component of $\boldsymbol{u}_1$, Eqs. (\ref{24}), (\ref{25}) and Eq. (\ref{30}) can be reduced to three equations for the perturbed variables namely the vertical
	component of the velocity $w_1$, and the concentration $n_1$:
	
	\begin{equation}\label{32}
	S_c^{-1}\nabla^2\left(\frac{\partial{w_1}}{\partial t}\right)=-\boldsymbol{\nabla}^4 w_1-Ra_b\nabla_h^2n_1+Ra_TLe\nabla_h^2T_1,
	\end{equation}
	\begin{equation}\label{33}
	\frac{\partial n_1}{\partial t}+w_1\frac{dn_s}{dz}+\Gamma_1(z)\int_1^zn_1 dz'+\Gamma_2(z)n_1+\Gamma_3(z)\frac{\partial n_1}{\partial z}=\nabla^2 n_1, 
	\end{equation}
	\begin{equation}\label{34}
	\frac{\partial{T_1}}{\partial t}+w_1\frac{dT_S}{dz}=Le\boldsymbol{\nabla}^2T_1,
	\end{equation}      
	
	These quantities can then be decomposed into normal modes such that      
	
	\begin{equation}\label{35}
	\begin{pmatrix}
	w_1\\n_1\\T_1
	\end{pmatrix}
	=
	\begin{pmatrix}
	W(z)\\\Theta(z)\\T(z)
	\end{pmatrix}
	+\exp{[\gamma t+i(lx+my)]},
	\end{equation} 
	
	where the amplitudes $W (z)$, $\Theta(z)$ and $T(z)$ describe the unknown variation concerning $z$ of the vertical velocity, the vertical vorticity and the concentration respectively, $l$ and $m$ are the dimensionless wave numbers in the $x$ and $y$ directions respectively,
	finally, $\gamma$ is the complex growth rate of the disturbances.     
	
	The linear stability equations become
	\begin{align}\label{36}
	\nonumber\left(\gamma Sc^{-1}+k^2-\frac{d^2}{dz^2}\right)\left( \frac{d^2}{dz^2}-k^2\right)W(z)=Ra_bk^2\Theta(z)\\-Ra_TLek^2T(z),
	\end{align}
	\begin{align}\label{37}
	\nonumber\Gamma_1(z)\int_1^z\Theta(z') dz'+\left(\sigma+k^2+\Gamma_2(z)+\Gamma_3(z)\frac{d}{dz}-\frac{d^2}{dz^2}\right)\Theta(z)\\=-\frac{dn_s}{dz}W(z), 
	\end{align}
	\begin{equation}\label{38}
	Le\boldsymbol{\nabla}^2T(z)-\gamma T(z)=\frac{dT_s}{dz}W(z),
	\end{equation}  
	
	subject to the boundary conditions
	\begin{align}\label{39}
	\nonumber W(z)=\frac{dW(z)}{dz}=\Gamma_2(z)\int_1^z\Theta(z')dz'+2\Gamma_3(z)\Theta(z)\\-2\frac{d\Theta(z)}{dz}=0,~~~at~~z=0,1,
	\end{align}
	\begin{equation}\label{40}
	T(z)=0,~~~at~~z=0,1.
	\end{equation}
	
	Here, $k=\sqrt{(l^2+m^2)}$ is the overall nondimensional wavenumber. Equations (\ref{36})-(\ref{40}) form an eigen value problem for $\gamma$ as a functions of the dimensionless parameters $l, m, V_c , \kappa$ and $R$. The basic state becomes unstable whenever $Re(\gamma) > 0$.
	Introducing a new variable
	
	\begin{equation}\label{41}
	\Phi(z)=\int_1^z\Theta(z')dz',
	\end{equation}
	the linear stability equations (after writing $D= d/dz$) become
	\begin{align}\label{42}
	\nonumber\left(\gamma S_c^{-1}+k^2-D^2\right)\left( D^2-k^2\right)W(z)=Ra_bk^2D\Phi(z)\\-Ra_tLek^2T(z),
	\end{align}
	
	\begin{align}\label{43}
	\nonumber\Gamma_1(z)\Phi(z)+(\sigma+k^2+\Gamma_2(z))D\Phi(z)+\Gamma_3(z)D^2\Phi(z)\\-D^3\Phi(z)=-Dn_sW(z), 
	\end{align}
	\begin{equation}\label{44}
	Le(D^2-k^2)T(z)=\gamma T(z)+DT_sW(z),
	\end{equation} 
	
	with the boundary conditions 
	
	\begin{align}\label{45}
	\nonumber W(z)=DW(z)=\Gamma_2(z)\Phi(z)+2\Gamma_3(z)D\Phi(z)\\-2D^2\Phi(z)=0,~~~at~~z=0,1,
	\end{align}
	\begin{equation}\label{46}
	T(z)=0,~~~at~~z=0,1,
	\end{equation}
	
	and extra boundary condition is
	\begin{equation}\label{47}
	\Phi(z)=0~~~at~~z=1,
	\end{equation}
	which is follows from Eq. (\ref{41}).
	
	The Eqs.~(\ref{42})-(\ref{44}) along with the boundary conditions can be written in the following matrix form
	
	\begin{equation}\label{48}
	A(k)Y=\gamma B(k)Y,
	\end{equation}	
	where $Y=(W,\Phi,T)$. $A(k)$ and $B(k)$ are two linear differential operators which are dependent on the control parameters $Ra_b$, $Ra_T$, $S_c$, $L_e$, $V_c$ and $\kappa$.
	
	%%%%%%%%%%%%%%%%%%%%%%%%%%%%%%%%%%%%%%%%%%%%%%%%%%%%%%%%%%%%%%%%%%%%%%%%%%%
	
	\section{SOLUTION technique and NUMERICAL RESULTS}
	To analyze the linear stability, Eq.~(\ref{48}) is solved by the Newton-Raphson-Kantorovich (NRK) iterative scheme (see Cash et al.~\cite{19cash1980}). This numerical scheme allows us to calculate the growth rate, Re$(\gamma)$, and neutral stability curves in the $(k, Ra)$-plane for a specific set of parameters where $Ra=Ra_b/Ra_T$. The neutral curve, denoted as $Ra^{(n)}(k)$, where $n$ is an integer greater than or equal to 1, consists of an infinite number of branches. Each branch represents a possible solution to the linear stability problem for the given parameter set. Among these branches, the one with the lowest value of $Ra$ is considered the most significant, and the corresponding bioconvective solution is identified as $(k_c, Ra_c)$. This particular solution is referred to as the most unstable solution. By utilizing the equation $\lambda_c=2\pi/k_c$, where $\lambda_c$ represents the wavelength of the initial disturbance, we can determine the wavelength associated with the most unstable solution. This wavelength provides valuable information about the characteristic pattern size of the bioconvection phenomenon.
	
	In this study, we recognize the complexity of exploring the entire parameter space due to the wide range of values that each parameter can take. To systematically investigate the impact of phototactic microorganisms on thermal convection, we decided to keep certain parameters constant while varying others. This approach allows us to focus on specific aspects of the system and gain a deeper understanding of their individual influence on the onset of bioconvection.
	
	We fix $S_c=20$, $G_t=0.8$, $L_e=1$ and $V_c=10$ as constant parameter values. These choices are made to maintain consistency and isolate the effects of other parameters. The parameters related to the absorption coefficient $\kappa$ and cell swimming speed $V_c$ are varied to observe their specific impacts. We consider two values for $\kappa$, namely 0.5 and 1.0, representing different light absorption characteristics of the microorganisms.
	
	Here, we determined the bioconvective Rayleigh number, $Ra_b$, at the onset of bioconvection as a function of wavenumber, $k$, for different thermal Rayleigh numbers, $Ra_T$.

	\begin{figure}[!htbp]
		\includegraphics[width=8cm]{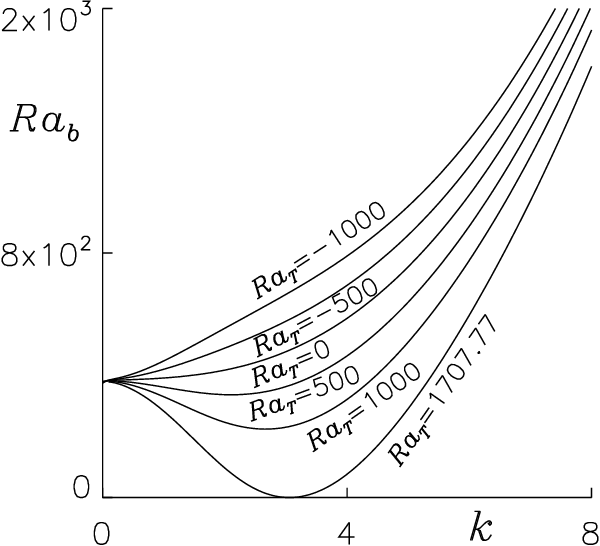}
		\caption{\label{fig2} The marginal stability curves $Ra_b$ vs $k$ for different levels of $Ra_T$ when $V_c=10, \kappa=0.5$, and $G_c=0.8$.}
	\end{figure}

	\begin{figure}[!htbp]
		\includegraphics[width=8cm]{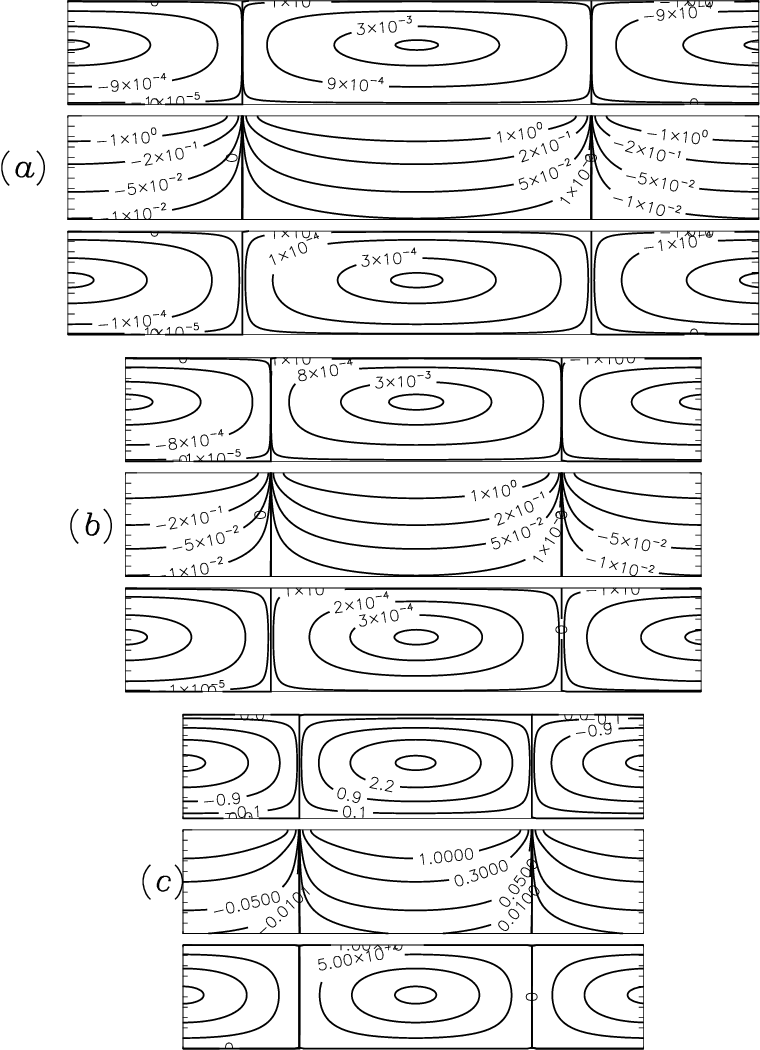}
		\caption{\label{fig3} Combined ﬂow patterns obtained via the perturbed velocity (top) and concentration (middle) and isotherms (bottom) component. The depicted patterns are shown in (a), (b) and (c) are obtained for the case of $Ra_T=-1000$, $0$ and $1000$ with $V_c=10, \kappa=0.5$, and $G_c=0.8$.}
	\end{figure}
	
	%%%%%%%%%%%%%%%%%%%%%%%%%%%%%%%%%%%%%%%%%%%%%%%%%%%%%%%%%%%%%%%%%%%%%%%%%%%%
	
	\subsection{$\kappa=0.5$}
	
	\begin{figure}[!ht]
		\includegraphics[width=8cm]{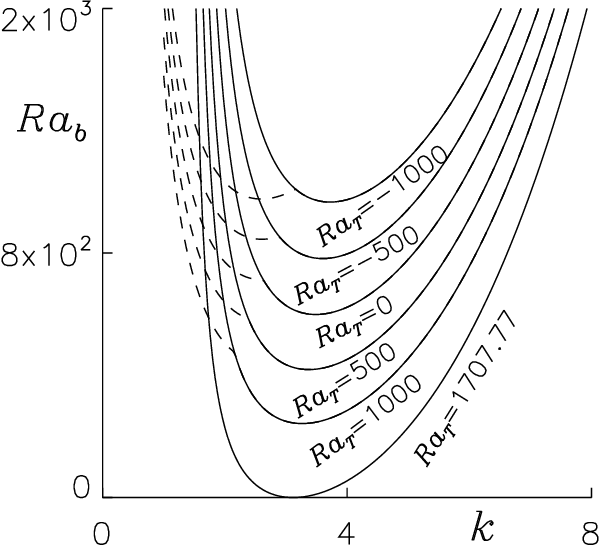}
		\caption{\label{fig4} The marginal stability curves $Ra_b$ vs $k$ for different $Ra_T$ when $V_c=10, \kappa=0.5$, and $G_c=0.66$.}
	\end{figure}
	
	\begin{figure}[!ht]
		\includegraphics[width=8cm]{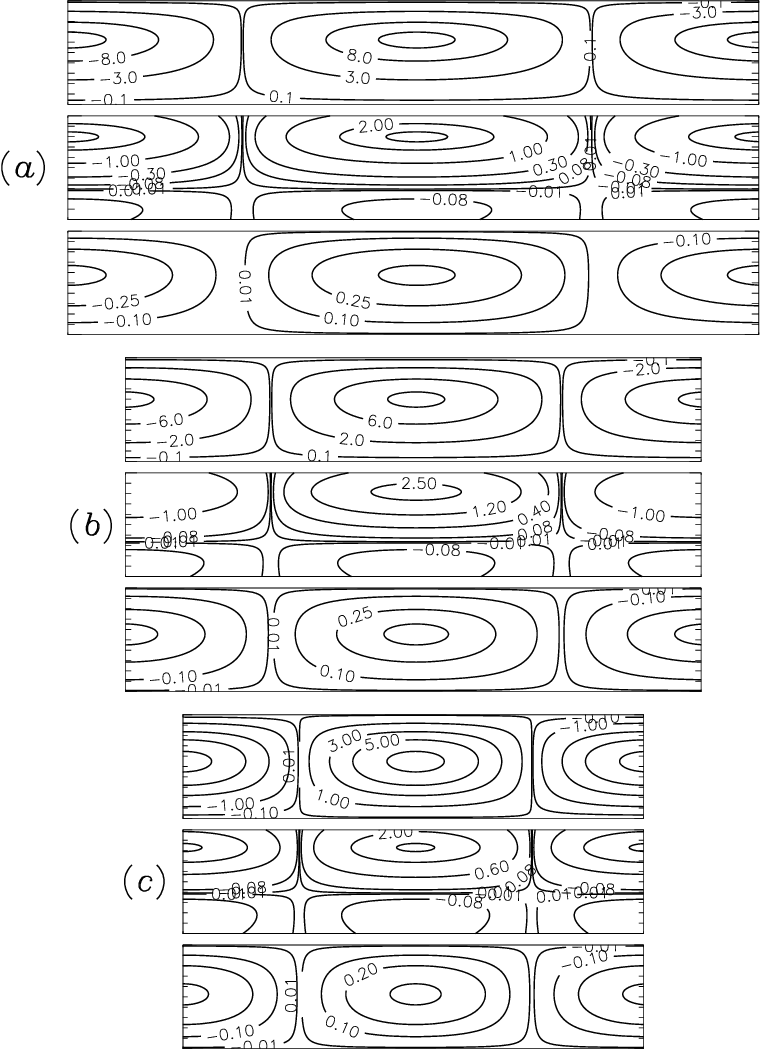}
		\caption{\label{fig5}  Combined ﬂow patterns obtained via the perturbed velocity (top) and concentration (middle) and isotherms (bottom) component. The depicted patterns are shown in (a), (b) and (c) are obtained for the case of $Ra_T=-1000$, $0$ and $1000$ with $V_c=10, \kappa=0.5$, and $G_c=0.66$.}
	\end{figure}
	
	%%%%%%%%%%%%%%%%%%%%%%%%%%%%%%%%%%%%%%%%%%%%%%%%%%%%%%%%%%%%%%%%%%%%%%%%%%%
	
	\begin{figure}[!ht]
		\includegraphics[width=8cm]{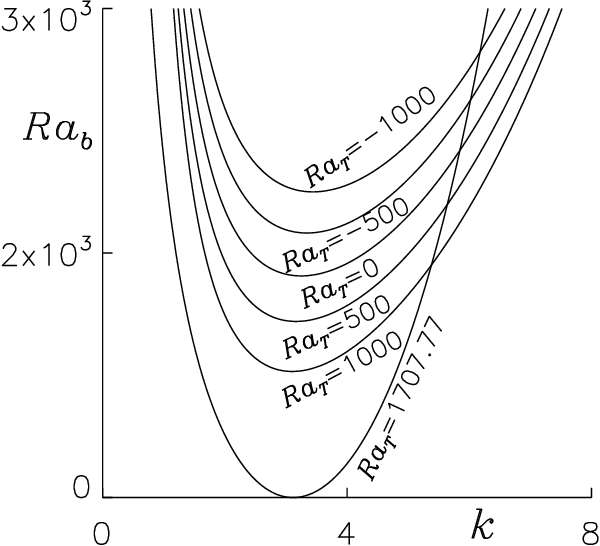}
		\caption{\label{fig6} The marginal stability curves $Ra_b$ vs $k$ for different $Ra_T$ when $V_c=10, \kappa=0.5$, and $G_c=0.63$.}
	\end{figure}

	In Fig.~\ref{fig2}, we present the marginal stability curves with bioconvective Rayleigh number $Ra_b$ vs wavenumber $k$ for $V_c=10$, $\kappa=0.5$, and $G_c=0.8$, while varying the thermal Rayleigh number ($Ra_T$) at different levels. For $Ra_T=1707.77$, the bioconvective Rayleigh number $Ra_b=0$ with pattern wavelength $\lambda=2.01$. As the $Ra_T$ decreases the critical bioconvective Rayleigh number and wavelength increase. At $Ra_T=0$, the critical wavelength becomes infinite and this behaviour continues for the rest of the thermal Rayleigh numbers $Ra_T$.

	For $G_c=0.66$, the marginal stability curves for different thermal Rayleigh numbers are shown in Fig.~\ref{fig4}, where the other governing parameters $V_c=10$, $\kappa=0.5$ are held fixed. Here, the zero critical bioconvective Rayleigh number occurs at $Ra_T=1707.77$. As the thermal Rayleigh number decreases up to $Ra_T=-1000$, the critical bioconvective Rayleigh number increases and the corresponding pattern wavelength decreases. Here, the oscillatory branch of the marginal stability curve bifurcates from the stationary branch but the most unstable solution occurs at the stationary branch which implies the stationary solution except $Ra_T=1707.77$.

	In Fig.~\ref{fig5}, the marginal stability curves for $G_c=0.63$ with different levels of thermal Rayleigh numbers are shown, where $V_c=10$ and $\kappa=0.5$ are kept hold. Here, the critical bioconvective Rayleigh number increases and the wavelength decreases as the thermal Rayleigh number decreases from $Ra_T=1707.77$ to -1000. In this section, stationary solutions are observed for different levels of critical light intensity.
	
	%\begin{table}[h]
	%	\caption{\label{tab3}The numerical values of the bioconvective solutions for $V_c=10$ and $V_c=15$, with an increase in $B$, while keeping other parameters constant.}
	%	\begin{ruledtabular}
	%		\begin{tabular}{ccccccccc}
	%			$G_c$ & $V_c$ & $\kappa$ & $Ra_T$ & $\lambda^c$ & $Ra_{b}^c$ & $Im(\sigma)$ \\
	%			\hline
	%			\vspace{-0.2cm}\\
	%			
	%			0.8 & 10 & 0.5 & 1707.77 & 2.01 & 0 & 0  \\
	%			0.8 & 10 & 0.5 & 1000 & 2.72 & 224.04 & 0 \\
	%			0.8 & 10 & 0.5 & 500 & 2.91 & 336.36 & 0 \\
	%			0.8 & 10 & 0.5 & 0 & $\infty$ & 374.76 & 0 \\
	%			0.8 & 10 & 0.5 & -500 & $\infty$ & 375.04 & 0 \\
	%			0.8 & 10 & 0.5 & -1000 & $\infty$ & 375.32 & 0 \\
	%			
	%			0.66 & 10 & 0.5 & 1707.77 & 2.01 & 0 & 0  \\
	%			0.66 & 10 & 0.5 & 1000 & 1.93 & 242.16 & 0 \\
	%			0.66 & 10 & 0.5 & 500 & 1.86 & 418.79 & 0 \\
	%			0.66 & 10 & 0.5 & 0 & 1.81 & 599.02 & 0 \\
	%			0.66 & 10 & 0.5 & -500 & 1.74 & 781.99 & 0 \\
	%			0.66 & 10 & 0.5 & -1000 & 1.69 & 966.84 & 0 \\
	%			
	%			0.63 & 10 & 0.5 & 1707.77 & 2.01 & 0 & 0  \\
	%			0.63 & 10 & 0.5 & 1000 & 2.01 & 773.62 & 0 \\
	%			0.63 & 10 & 0.5 & 500 & 1.99 & 1080.50 & 0 \\
	%			0.63 & 10 & 0.5 & 0 & 1.94 & 1359.80 & 0 \\
	%			0.63 & 10 & 0.5 & -500 & 1.87 & 1623.15 & 0 \\
	%			0.63 & 10 & 0.5 & -1000 & 1.83 & 1875.40 & 0\\
	%			
	%		\end{tabular}
	%	\end{ruledtabular}
	%	\footnotetext[1]{The result indicates the existence of the oscillatory branch of the marginal stability curve.}
	%	\footnotetext[2]{The result indicates that the most unstable solution occurs on the oscillatory branch.}
	%\end{table}

	%%%%%%%%%%%%%%%%%%%%%%%%%%%%%%%%%%%%%%%%%%%%%%%%%%%%%%%%%%%%%%%%%%%%%%%%%%
	
	\subsection{$\kappa=1$}

	\begin{figure}[!ht]
		\includegraphics[height=8cm]{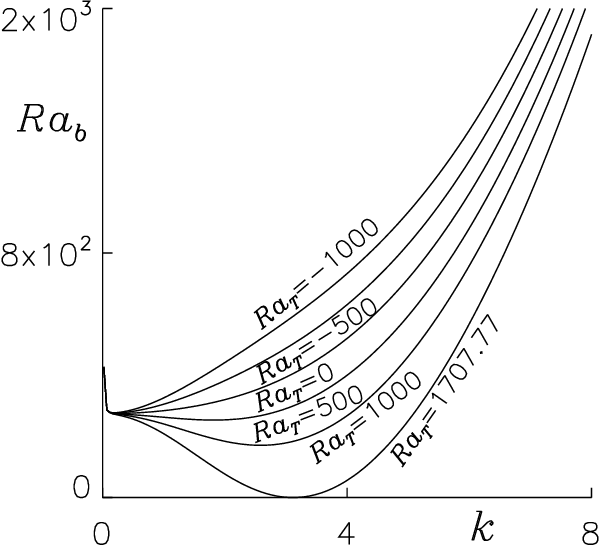}
		\caption{\label{fig7} The marginal stability curves $Ra_b$ vs $k$ for different $Ra_T$ when $V_c=10, \kappa=1$, and $G_c=0.8$.}
	\end{figure}
	
	\begin{figure}[!ht]
		\includegraphics[width=8cm]{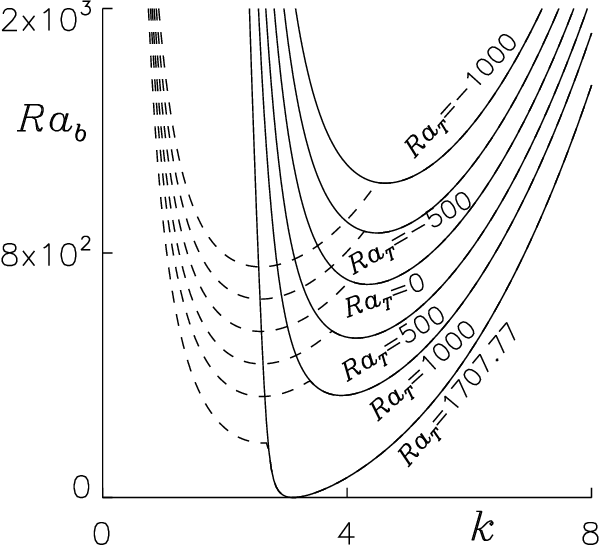}
		\caption{\label{fig8} The marginal stability curves $Ra_b$ vs $k$ for different $Ra_T$ when $V_c=10, \kappa=1$, and $G_c=0.53$.}
	\end{figure}

	\begin{figure}[!ht]
		\includegraphics[width=8cm]{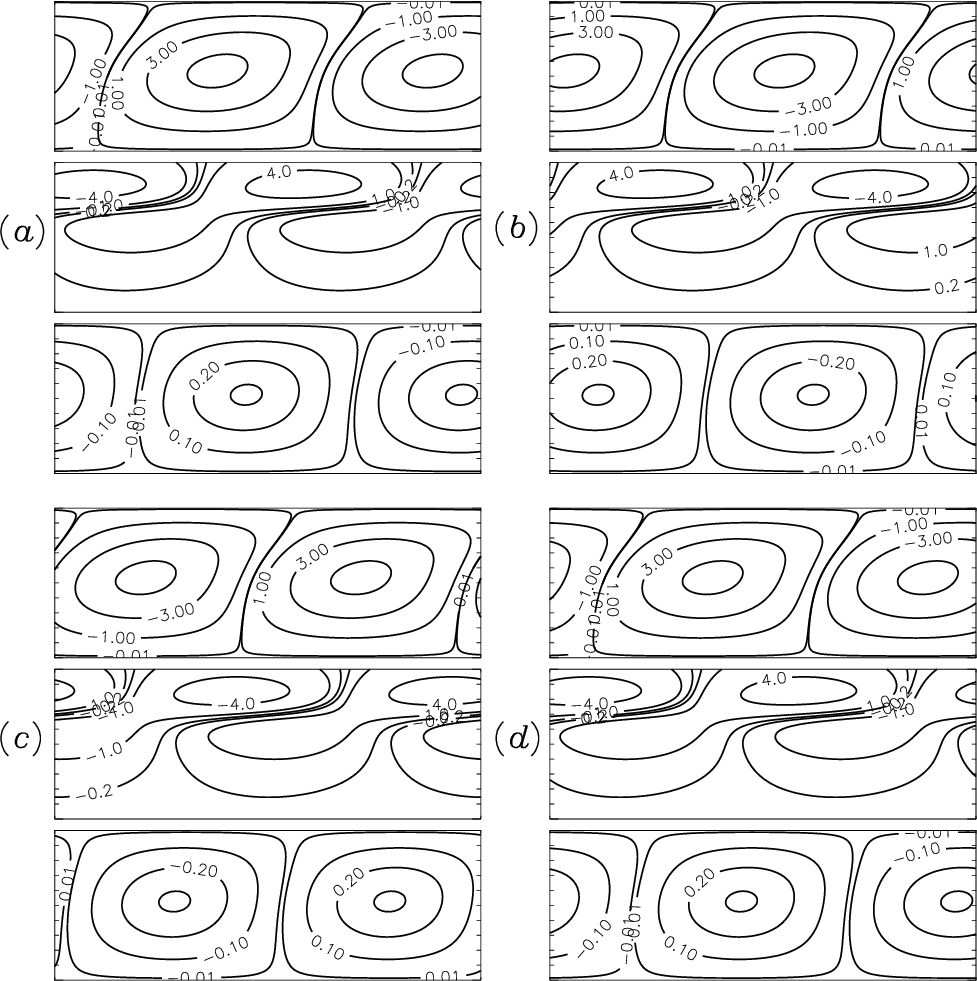}
		\caption{\label{fig9} Combined ﬂow patterns obtained via the perturbed velocity (top) and concentration (middle) and isotherms (bottom) component, during a cycle of oscillation at the onset of bioconvective overstability. The depicted patterns shown in (a)–(d) are obtained for the case of $Ra_T=1000$ with $V_c=10, \kappa=1$, and $G_c=0.53$, and the time intervals in between them are equal.}
	\end{figure}
	\begin{figure}[!ht]
		\includegraphics[width=8cm]{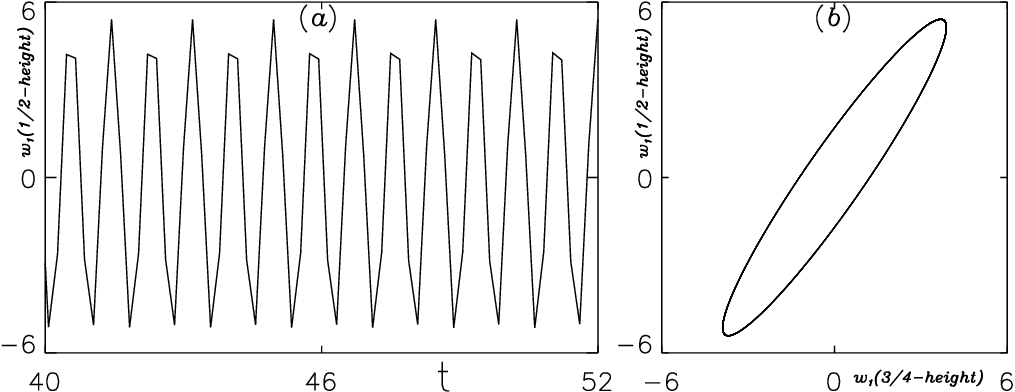}
		\caption{\label{fig10} (a) Time-evolving perturbed fluid velocity
			$w_1$ and (b) limit cycle corresponding to the overstable branch of the neutral curve for $Ra_T=1000$ shown in Fig.~\ref{fig8}. Fixed parameter values are cell swimming speed $V_c=10$, optical
			depth $\kappa= 1$, and critical total intensity $G_c=0.53$.}
	\end{figure}

	\begin{figure}[!ht]
		\includegraphics[width=8cm]{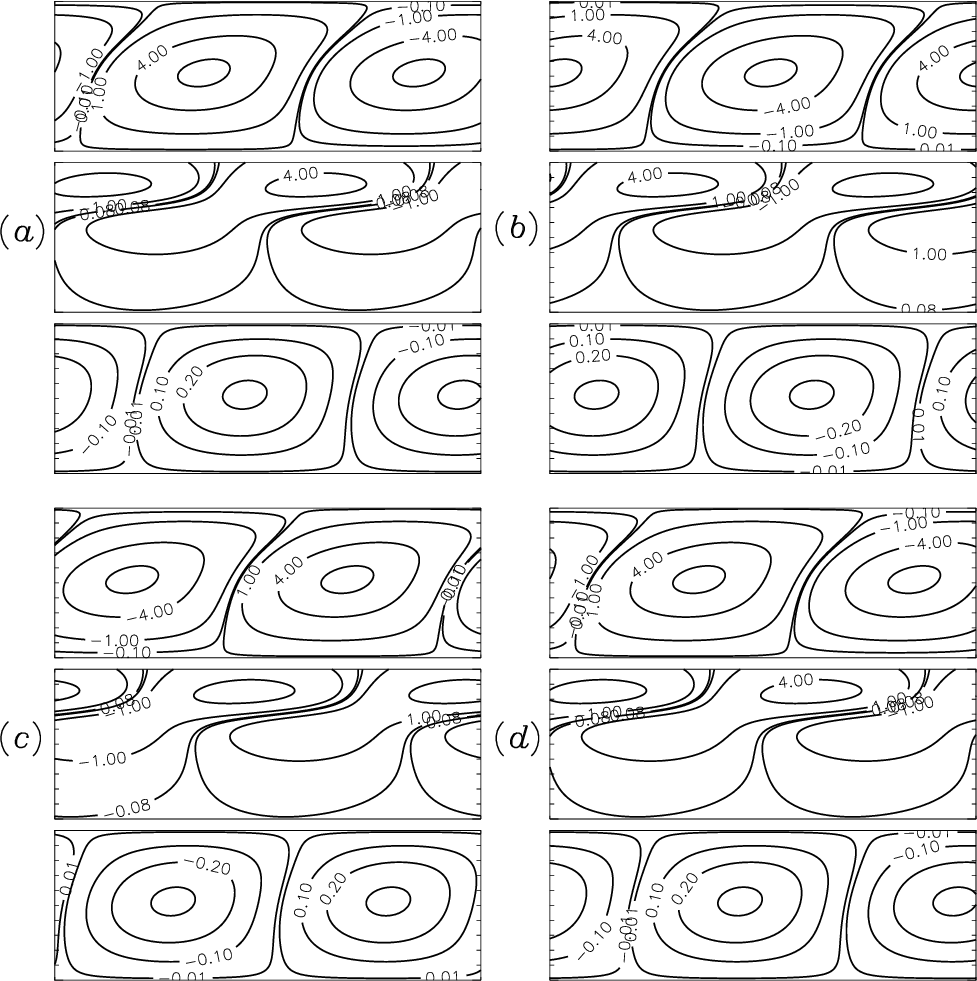}
		\caption{\label{fig11} Combined ﬂow patterns obtained via the perturbed velocity (top) and concentration (middle) and isotherms (bottom) component, during a cycle of oscillation at the onset of bioconvective overstability. The depicted patterns shown in (a)–(d) are obtained for the case of $Ra_T=0$ with $V_c=10, \kappa=1$, and $G_c=0.53$, and the time intervals in between them are equal.}
	\end{figure}

	\begin{figure}[!ht]
		\includegraphics[width=8cm]{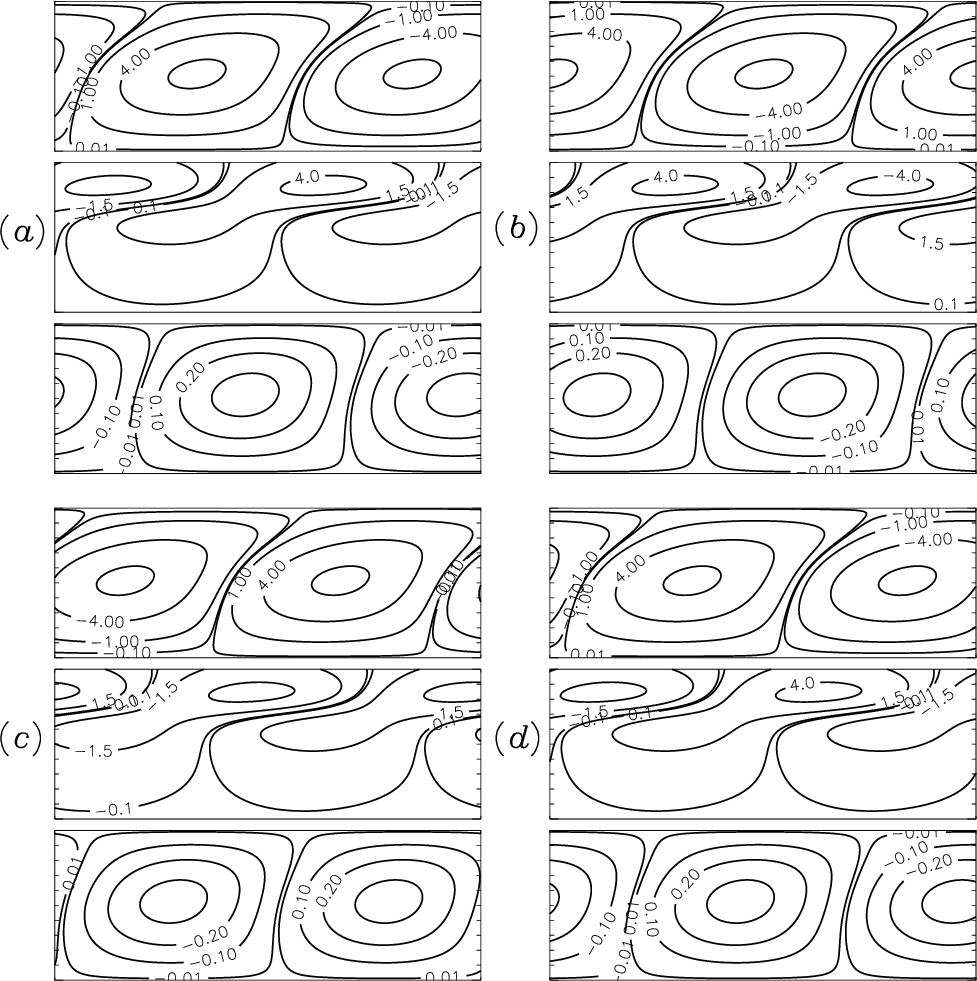}
		\caption{\label{fig12} Combined ﬂow patterns obtained via the perturbed velocity (top) and concentration (middle) and isotherms (bottom) component, during a cycle of oscillation at the onset of bioconvective overstability. The depicted patterns shown in (a)–(d) are obtained for the case of $Ra_T=-1000$ with $V_c=10, \kappa=1$, and $G_c=0.53$, and the time intervals in between them are equal.}
	\end{figure}
	
	\begin{figure}[!ht]
		\includegraphics[width=8cm]{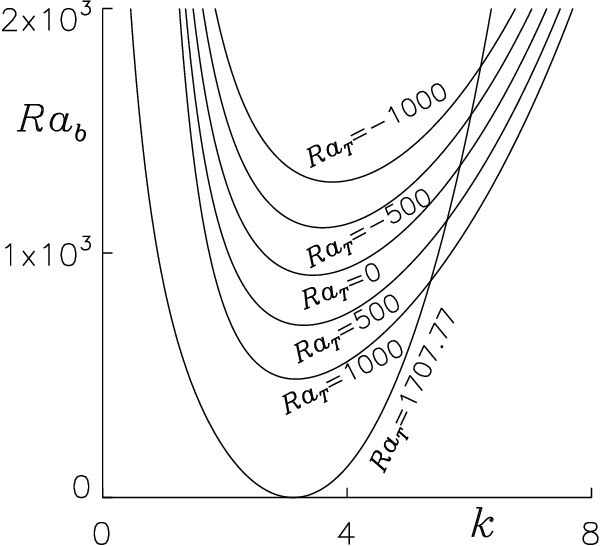}
		\caption{\label{fig13} The marginal stability curves $Ra_b$ vs $k$ for different $Ra_T$ when $V_c=10, \kappa=1$, and $G_c=0.595$.}
	\end{figure}
	
	In Fig.~\ref{fig11}, we present the marginal stability curves with bioconvective Rayleigh number $Ra_b$ vs wavenumber $k$ for $V_c=10$, $\kappa=1$, and $G_c=0.8$, while varying the thermal Rayleigh number ($Ra_T$) at different levels. For $Ra_T=1707.77$, the bioconvective Rayleigh number $Ra_b=0$ with pattern wavelength $\lambda=2.01$. As the $Ra_T$ decreases from 1707.77 to -1000, the critical bioconvective Rayleigh number and the pattern wavelength increase.
	
	For $G_c=0.53$, the marginal stability curves for different thermal Rayleigh numbers are shown in Fig.~\ref{fig8}, where the other governing parameters $V_c=10$, $\kappa=1$ are held fixed. Here, the zero critical bioconvective Rayleigh number occurs at $Ra_T=1707.77$ and an oscillatory branch bifurcates from the stationary branch of the marginal stability curve, but the most unstable solution occurs at the stationary branch which implies a stationary bioconvective solution. As the thermal Rayleigh number decreases up to $Ra_T=1000$, the critical bioconvective Rayleigh number and corresponding pattern wavelength increase. Here, also the oscillatory branch from the marginal stability curves bifurcate from the stationary branch, but the most unstable solution occurs at the oscillatory branch which implies the overstable solution. Thus, the onset of overstability is at $k_c=2.32$ and $Ra_{bc}=330.58$. At this point, two complex conjugate eigenvalues $Im(\gamma)=\pm 7.14i$ are found. The transition observed here is known as Hopf bifurcation. The bioconvective flow patterns, which correspond to the complex conjugate pair of eigenvalues, are mirror images of each other. The period of oscillation is $2\pi/Im(\gamma)=0.88$ units. The bioconvective fluid motions become fully nonlinear on a timescale substantially less than the predicted period of overstability. Thus, the convection cells and flow patterns during one cycle of oscillation can be visualized via the perturbed eigenmodes $w_1$, $n_1$ and $T_1$, respectively (see Fig.~\ref{fig9}). It reveals that a travelling wave solution is moving toward the left of the figure. Fig. \ref{fig10} shows the predicted time-evolving perturbed fluid velocity component $w_1$ [Fig. \ref{fig10}(a)] and its corresponding phase diagram [Fig. \ref{fig10}(b)] at $k_c=2.32$. It is again worth noting that the period of oscillation, that is, $2\pi/Im(\gamma)$ is the bifurcation (control) parameter, and hence, the bioconvective flow destabilization gives birth to a limit cycle (i.e., isolated cycle) [see Fig. \ref{fig10}(b)]. The birth of a limit cycle due to flow destabilization is again recognized as the Hopf bifurcation via bifurcation analysis. Based on the linear stability theory, the supercritical nature of this Hopf bifurcation finally leads it into a stable limit cycle. As $Ra_T$ is decreased up to -1000, the same trend continues and the oscillatory branch retains the most unstable mode resulting in the bioconvective solution being overstable.
	
	\begin{table}[h]
		\caption{\label{tab3}The numerical values of the bioconvective solutions for $V_c=10$ and $\kappa=0.5$, with different values of the critical light intensities, while keeping other parameters constant.}
		\begin{ruledtabular}
			\begin{tabular}{ccccccccc}
				$G_c$ & $V_c$ & $\kappa$ & $Ra_T$ & $\lambda^c$ & $Ra_{b}^c$ & $Im(\gamma)$ \\
				\hline
				\vspace{-0.2cm}\\
				
				0.8 & 10 & 0.5 & 1707.77 & 2.01 & 0 & 0  \\
				0.8 & 10 & 0.5 & 1000 & 2.72 & 224.04 & 0 \\
				0.8 & 10 & 0.5 & 500 & 2.91 & 336.36 & 0 \\
				0.8 & 10 & 0.5 & 0 & $\infty$ & 374.76 & 0 \\
				0.8 & 10 & 0.5 & -500 & $\infty$ & 375.04 & 0 \\
				0.8 & 10 & 0.5 & -1000 & $\infty$ & 375.32 & 0 \\
				
				0.66 & 10 & 0.5 & 1707.77 & 2.01 & 0 & 0  \\
				0.66 & 10 & 0.5 & 1000 & 1.93 & 242.16 & 0 \\
				0.66 & 10 & 0.5 & 500 & 1.86 & 418.79 & 0 \\
				0.66 & 10 & 0.5 & 0 & 1.81 & 599.02 & 0 \\
				0.66 & 10 & 0.5 & -500 & 1.74 & 781.99 & 0 \\
				0.66 & 10 & 0.5 & -1000 & 1.69 & 966.84 & 0 \\
				
				0.63 & 10 & 0.5 & 1707.77 & 2.01 & 0 & 0  \\
				0.63 & 10 & 0.5 & 1000 & 2.01 & 773.62 & 0 \\
				0.63 & 10 & 0.5 & 500 & 1.99 & 1080.50 & 0 \\
				0.63 & 10 & 0.5 & 0 & 1.94 & 1359.80 & 0 \\
				0.63 & 10 & 0.5 & -500 & 1.87 & 1623.15 & 0 \\
				0.63 & 10 & 0.5 & -1000 & 1.83 & 1875.40 & 0\\
				
			\end{tabular}
		\end{ruledtabular}
	\end{table}

	\begin{table}[h]
		\caption{\label{tab3}The numerical values of the bioconvective solutions for $V_c=10$ and $\kappa=1$, with different values of the critical light intensities, while keeping other parameters constant.}
		\begin{ruledtabular}
			\begin{tabular}{ccccccccc}
				$G_c$ & $V_c$ & $\kappa$ & $Ra_T$ & $\lambda^c$ & $Ra_{b}^c$ & $Im(\gamma)$ \\
				\hline
				\vspace{-0.2cm}\\
				
				0.8 & 10 & 1 & 1707.77 & 2.01 & 0 & 0  \\
				0.8 & 10 & 1 & 1000 & 2.42 & 171 & 0 \\
				0.8 & 10 & 1 & 500 & 3.49 & 253.56 & 0 \\
				0.8 & 10 & 1 & 0 & 23.27 & 275.28 & 0 \\
				0.8 & 10 & 1 & -500 & 28.56 & 276.60 & 0 \\
				0.8 & 10 & 1 & -1000 & 36.96 & 277.49 & 0 \\
				
				0.53 & 10 & 1 & 1707.77 & 2.01 & 0 & 0  \\
				0.53 & 10 & 1 & 1000 & 2.42 & 330.58 & 7.14 \\
				0.53 & 10 & 1 & 500 & 2.42 & 437.28 & 9.09 \\
				0.53 & 10 & 1 & 0 & 2.42 & 543.52 & 10.67 \\
				0.53 & 10 & 1 & -500 & 2.42 & 649.34 & 12.04 \\
				0.53 & 10 & 1 & -1000 & 2.45 & 754.82 & 13.21 \\
				
				0.5 & 10 & 1 & 1707.77 & 2.01 & 0 & 0  \\
				0.5 & 10 & 1 & 1000 & 1.99 & 484.57 & 0 \\
				0.5 & 10 & 1 & 500 & 1.92 & 704.63 & 0 \\
				0.5 & 10 & 1 & 0 & 1.83 & 909.26 & 0 \\
				0.5 & 10 & 1 & -500 & 1.75 & 1103.77 & 0 \\
				0.5 & 10 & 1 & -1000 & 1.67 & 1290.50 & 0\\
				
			\end{tabular}
		\end{ruledtabular}
	\end{table}

	\subsection{Effects of cell swimming speed}
	
\begin{figure}[!htbp]
		\includegraphics[width=8cm]{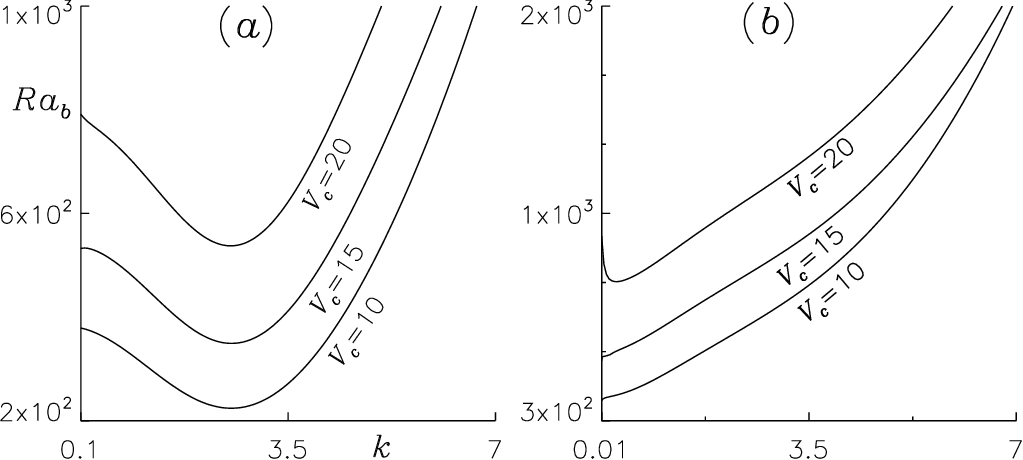}
		\caption{\label{fig14} The marginal stability curves $Ra_b$ vs $k$ for different $V_c$ when (a) $Ra_T=1000$, and (b) $Ra_T=-1000$. Here, the other governing parameters are $\kappa=0.5$ and $G_c=0.8$.}
	\end{figure}

	The effects of cell swimming speed at the onset of thermo-photo-bioconvection are also investigated. First, we discuss the case of $\kappa=0.5$. The Figs~\ref{fig14}(a) and \ref{fig14}(b) showcase the marginal stability curves for different cell swimming speeds at the thermal Rayleigh number $Ra_T=1000$ and $Ra_T=-1000$ when $G_c=0.8$. For both of the cases, the critical Bioconvective Rayleigh number increases as the cell swimming speed increases.
	
\begin{figure}[!htbp]
		\includegraphics[width=8cm]{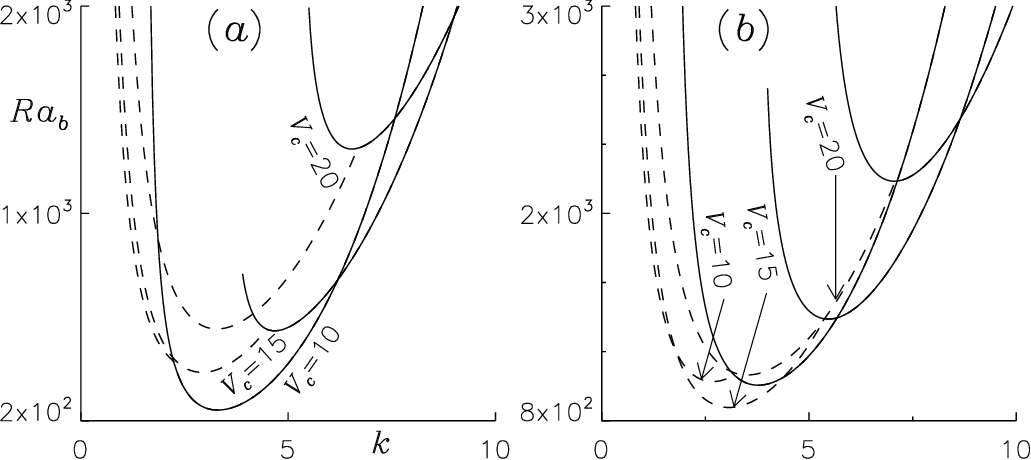}
		\caption{\label{fig15} The marginal stability curves $Ra_b$ vs $k$ for different $V_c$ when (a) $Ra_T=1000$, and (b) $Ra_T=-1000$. Here, the other governing parameters are $\kappa=0.5$ and $G_c=0.66$.}
	\end{figure}

	The marginal stability curves at $Ra_T=1000$ and $Ra_T=-1000$ are shown in Fig. \ref{fig15} for $V_c=10$ to 20. Here $G_c=0.66$ is a fixed constant value. For the case of $Ra_T=1000$ the critical bioconvective Rayleigh number increases as the cell swimming speed increases. Now come to the case of $Ra_T=-1000$. Here lower critical bioconvective Rayleigh number is observed for $V_c=15$ as compared to $V_c=10$ but $Ra_b^c$ increases when $V_c$ increases to 20. In both of the situations, the oscillatory solutions are also observed. 
	
\begin{figure}[!htbp]
		\includegraphics[width=8cm]{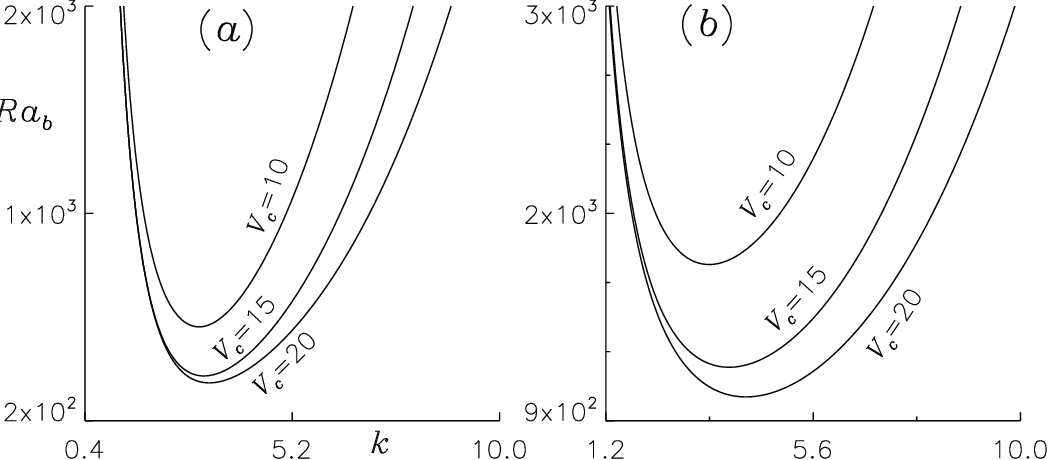}
		\caption{\label{fig16} The marginal stability curves $Ra_b$ vs $k$ for different $V_c$ when (a) $Ra_T=1000$, and (b) $Ra_T=-1000$. Here, the other governing parameters are $\kappa=0.5$ and $G_c=0.63$.}
	\end{figure}

	For $\kappa=0.5$ and $G_c=0.63$, the marginal stability curves of different cell swimming speeds are shown in Fig.~\ref{fig16}. Here Fig.~\ref{fig16}(a) is related to the $Ra_T=1000$ and Fig.~\ref{fig16}(b) is related to the $Ra_T=-1000$. The critical bioconvective Rayleigh number decreases as the cell swimming speed increases from 10 to 20 in both of the cases.
	
	For the case of $\kappa=1$, similar patterns occur.

	%%%%%%%%%%%%%%%%%%%%%%%%%%%%%%%%%%%%%%%%%%%%%%%%%%%%%%%%%%%%%%%%%%%%%%%%%	
	
	\section{Conclusion}
	
	We introduce, for the first time, a thermal phototactic bioconvection model, incorporating the thermal impact on the system. Our focus is on investigating the influence of heating or cooling on the onset of phototactic bioconvection. The model adopts rigid no-slip boundary conditions on both the top and bottom of the suspension and linear perturbation theory is employed for the linear instability analysis.
	
	The linear stability analysis encompasses both stationary and oscillatory solutions, with the latter being predominantly observed when the microorganism sublayer is situated at three-quarters of the suspension's height. The oscillatory nature of the solutions transitions into stability as the thermal critical Rayleigh number increases. Furthermore, an increase (decrease) in the critical bioconvective Rayleigh number is noted as the thermal critical Rayleigh number increases (decreases). This observed relationship suggests that heating from below destabilizes the suspension while cooling from below stabilizes it.
	
	On the other hand, the critical pattern wavelength experiences an increase as the critical thermal Rayleigh number increases when the microorganism sublayer is positioned at the top of the suspension. However, a decrease in the critical light intensity results in a decrease in the critical pattern wavelength as the critical thermal Rayleigh number increases.

	%%%%%%%%%%%%%%%%%%%%%%%%%%%%-OM NAMAH SHIVAY-%%%%%%%%%%%%%%%%%%%%%%%%%%%%	
	
	\section*{ Availability of Data}
	The supporting data of this article is available within the article. 
	\nocite{*}
	\section*{REFERENCES}
	\bibliography{Thermo_bioconvection}
	
\end{document}